\documentclass[letterpaper, 10 pt, conference]{ieeeconf}  

\IEEEoverridecommandlockouts                              
\usepackage{cite}
\usepackage{amsmath,amssymb,amsfonts}
\usepackage{graphicx}
\usepackage{textcomp}
\usepackage{xcolor}

\usepackage{epstopdf}
\usepackage{algorithm}
\usepackage{algorithmic}
\usepackage{booktabs}
\usepackage{multirow}

\DeclareMathOperator*{\argmin}{argmin}
\newtheorem{theorem}{\indent Theorem}
\newtheorem{remark}{\indent Remark}
\title{Distributed Noise Covariance Matrices Estimation in Sensor Networks}

\author{Jiahong Li, Nan Ma, and Fang Deng
\thanks{}
\thanks{Jiahong Li and Nan Ma are now with the School of Robotics, Beijing Union University, Beijing 100101, China. {\tt\small jqrjiahong@buu.edu.cn, xxtmanan@buu.edu.cn}}%
\thanks{Fang Deng is with the School of Automation, Beijing Institute of Technology, Beijing 100081, China. {\tt\small dengfang@bit.edu.cn}.}%
}

\begin{document}

\maketitle
\thispagestyle{empty}
\pagestyle{empty}

\begin{abstract}
Adaptive algorithms based on in-network processing over networks are useful for online parameter estimation of historical data (e.g., noise covariance) in predictive control and machine learning areas.
This paper focuses on the distributed noise covariance matrices estimation problem for multi-sensor linear time-invariant (LTI) systems.
Conventional noise covariance estimation approaches, e.g., auto-covariance least squares (ALS) method, suffers from the lack of the sensor's historical measurements and thus produces high variance of the ALS estimate.
To solve the problem, we propose the distributed auto-covariance least squares (D-ALS) algorithm based on the batch covariance intersection (BCI) method by enlarging the innovations from the neighbors.
The accuracy analysis of D-ALS algorithm is given to show the decrease of the variance of the D-ALS estimate.
The numerical results of cooperative target tracking tasks in static and mobile sensor networks are demonstrated to show the feasibility and superiority of the proposed D-ALS algorithm.
\end{abstract}

\section{Introduction}\label{sec:intro}
Recent advances in machine learning and information fusion have led to the formulation of increasingly demanding distributed estimation and inference problems, as discussed in \cite{He2019}.
The distributed estimation fusion methods in \cite{Cattivelli2010}, and especially the batch covariance intersection (BCI) approach (see \cite{Chong2001} and references therein) provided an upper bound on estimation accuracy without assuming any knowledge on the correlation between the estimates of sensors.
\cite{Sun2016} proposed a average consensus estimation algorithm based on a new BCI strategy.
However, the fusion methods above lack the consideration of the exact knowledge of the noise statistics, which is not plausible due to the mismatch of the nominal system or invalidity of offline calibration in many practical systems, e.g., low-cost integrated GPS/INS positioning systems \cite{TerryMoore2003}, energy-based source localization \cite{Deng2017} and fault tolerant systems \cite{Gertler1988}.

One effective approach is to use the historical open-loop data, which can be divided into several categories, e.g., correlation techniques \cite{Mehra1970,Belanger1974,Odelson2006,Akesson2008,Rajamani2009,Dunik2009,Abdel2010,Lima2013,Dunik2016}, Bayesian \cite{Sarkka2009,Li2019}, maximum likelihood \cite{Kashyap1970,Chen2017}, covariance matching \cite{Myers1976}, methods based on the minimax approach \cite{Verdu1984}, subspace methods \cite{Overschee1996} and prediction error methods \cite{Carew1974}.
An alternative approach that directly estimates the gain of a linear estimator has been developed in \cite{Carew1974,Wiberg2000,Dunik2009,Lima2013}.
The connections between two approaches were discussed in \cite{Dunik2016}.
The Bayesian and maximum likelihood methods are well suited to multi-model approaches, but are costly in terms of computation.
Covariance matching is a technique to provide biased estimates of the true covariances based on the residuals of the state estimates.
The minimax approach provides a fixed system whose worst performance among an assumed possible uncertainty set is the best possible.
The advantages and disadvantages of the approach have been discussed in \cite{Basar1972,Tugnait1979}.
The subspace methods formulate the estimation problem as projections of Hankel matrices and the model can be retrieved from the row and column spaces of the projected data matrix.
The prediction error methods reduces the parameter identification problem to the minimization of empirical average losses.
Among all the methods, the correlation methods can provide unbiased estimates with acceptable computational requirements even for high-dimensional systems \cite{Dunik2016}.
The correlation methods were firstly proposed by Mehra and B\'{e}langer in \cite{Mehra1970} and \cite{Belanger1974} as a three-step procedure, and were reformulated to a single-step procedure called the auto-covariance least-squares (ALS) method in \cite{Odelson2006}.
In the ALS method, the correlations between routine operating data formed a least-squares problem of the noise covariance matrix, whose solution was guaranteed by solving the semi-definite programming (SDP) problem.
The necessary and sufficient conditions for the uniqueness of the variance estimates for dependent state and measurement noise were presented in \cite{Akesson2008}.
The ALS problem with the estimation of a state noise disturbance structure was formulated in \cite{Rajamani2009}.
The optimal weight was formulated in the least-squares objective to ensure minimum variance in \cite{Dunik2016}.
However, the performance of the correlation methods would become poor if the time window size of open-loop measurements is small.

In this paper, the distributed noise covariance estimation problem over networks is formulated, and the distributed auto-covariance least squares (D-ALS) algorithm are proposed based on batch covariance intersection (BCI) method.
The estimation accuracy of the proposed algorithm can increase by fusing the innovations from the neighboring agents.
The theoretical analysis of the algorithm is also provided to shown the efficiency.
The simulation results of cooperative target tracking case show the superiority of the ALS-BCI algorithm in terms of the mean square error criterion.
\section{Preliminaries}
\label{sec:problem}
We consider a connected sensor network of $M$ agents modeled as an undirected graph $\mathcal {G}(\mathcal {V},\mathcal {E})$,
where the vertices set $\mathcal {V}={1,\ldots,M}$ corresponds to the agents and the edge set $\mathcal {E} \subset \mathcal {V} \times \mathcal {V}$ represents the communication links between the pairs of agents.
Agent $i$ can communicate with its neighbors
whose indexes are in the set $\mathcal{N}_{i} = \left\{j\in \mathcal {V}:(i,j)\in \mathcal {E}, i\ne j\right\}$ with cardinality $M_{i}=\|\mathcal{N}_{i}\|$.

In the sensor network, each agent observes the linear discrete time-invariant dynamic system $x_{k+1} = Fx_{k}+w_{k}$ with linear time-invariant measurement model $z_{i,k} = H_{i}x_{k} + v_{i,k}$.
where the vector $x_{k}\in \mathbb{R}^{n_{x}}$ and $z_{i,k}\in \mathbb{R}^{n_{z}}$ represent the state and the measurement of the $i^{th}$ agent at time instant $k\in \mathbb{N}^{+}$.
$F\in \mathbb{R}^{n_{x}\times n_{x}}$ and $H_{i}\in \mathbb{R}^{n_{z}\times n_{x}}$ are state-transitional and measurement-transitional matrix.
The variables $w_{k}$ and $v_{i,k}$ represent the process noise and measurement noise respectively, and are mutually independent following the zero-mean Gaussian statistics with probability $w_{k}\sim \mathcal{N}(0_{n_{x}\times 1},Q)$ and $v_{i,k}\sim \mathcal{N}(0_{n_{z}\times 1},R_{i})$ with unknown covariance matrices $Q\in \mathbb{R}^{n_{x}\times n_{x}}$ and $R_{i}\in \mathbb{R}^{n_{z}\times n_{z}}$.
$I_{n}$ and $0_{n}$ denote the identity matrix and the zeros matrix of dimension $n$ respectively.

According to the conventional distributed linear filtering algorithm, each agent updates local state estimate $\hat{x}_{i,k} = F\hat{x}_{i,k-1} + K_{i}e_{i,k}$ and state covariance estimate $P_{i,k} = (I-K_{i}H_{i})(FP_{i,k-1}^{-1}F^\mathrm{T}+Q)$ with estimation gain $K_{i}\in\mathbb{R}^{n_{x}\times n_{z}}$, and transmits them to its neighbors to fuse the global ones.
where $e_{i,k}= z_{i,k}-H_{i}F\hat{x}_{i,k-1}$ denote the innovation.
$K_{i}$ is designed as {K}alman gain $K_{i} =P_{i,k-1}H^\mathrm{T}_{i}(H_{i}(FP_{i,k-1}^{-1}F^\mathrm{T}+Q)H^\mathrm{T}_{i}+R_{i})^{-1}$ in terms of minimum mean square error (MMSE).
Noting that the noise covariance matrices $Q$ and $R_{i}$ are unknown and estimated by the auto-covariance least-squares method below.

Denote the residuals and residuals covariance of the $i^{th}$ agent as $\varepsilon_{i,k}=x_{i,k}-F\hat{x}_{i,k-1}$ and $P_{\varepsilon,i,k} = E[\varepsilon_{i,k}\varepsilon_{i,k}^\mathrm{T}]$ respectively, then the estimator of the residuals is deduced as
\begin{equation}\label{eps_i}
\varepsilon_{i,k} = \underbrace{(F-K_{i}H_{i}F)}_{\bar{F}_{i}}\varepsilon_{i,k-1} + \underbrace{[I_{n_{x}}-K_{i}H_{i},-K_{i}]}_{G_{i}}\underbrace{\left[
                                  \begin{array}{c}
                                    w_{k} \\
                                    v_{i,k} \\
                                  \end{array}
                                \right]}_{\bar{w}_{i,k}}
\end{equation}
\begin{equation}\label{Peps_i}
  P_{\varepsilon,i,k} = \bar{F_{i}}P_{\varepsilon,i,k-1}\bar{F_{i}}^\mathrm{T} + G_{i}\Sigma_{i} G_{i}^\mathrm{T}
\end{equation}
where $\Sigma = E(\bar{w}_{i,k}\bar{w}_{i,k}^\mathrm{T})=\left[
         \begin{array}{cc}
           Q & 0_{n_x\times n_z} \\
           0_{n_z\times n_x} & R_{i} \\
         \end{array}
       \right]$.
According to the Lyapunov equation $P_{\varepsilon} = \bar{F}P_{\varepsilon}\bar{F}^\mathrm{T} + G\Sigma G^\mathrm{T}$ in \cite{Simon2006}, the steady-state residual covariance solution exists if $\bar{F}$ is stable.
To ensure $\bar{F}$ is stable, the residual covariance $P_{\varepsilon}$ should satisfy $(P_{\varepsilon})_{s}=\big((I-\bar{F}\otimes\bar{F})^{-1}G\otimes G \big)\Sigma_{s}$ through the vectorization, where $\otimes$ denotes the Kronecker product, $A_s$ denotes the columnwise stacking of the matrix $A$ into a vector.

The innovations $e_{i,k}$ is deduced as
$e_{i,k} =  H_{i}\varepsilon_{i,k} + v_{i,k}$.
Then denote the auto-covariance $C^{i}_{e,0} = \mathbb{E}[e_{i,k}e_{i,k}^\mathrm{T}]$ and $C^{i}_{e,l} = \mathbb{E}[e_{i,k+l}e_{i,k}^\mathrm{T}]$ of the $i^{th}$ agent's innovation as
\begin{equation}\label{AC_inno_i}
  C^{i}_{e,l} = H_{i}\bar{F}^{l}P_{\varepsilon}H_{i}^\mathrm{T}- H_{i}\bar{F}^{l-1}FK_{i}R_{i} \quad l=1,2,\ldots,N-1
\end{equation}
where $C^{i}_{e,0} = H_{i}P_{\varepsilon,i}H_{i}^\mathrm{T}+ R_{i}$, $N$ is a user-defined parameter defining the maximum time-window lag.
It can be derived as
\begin{equation}\label{Ab}
  \mathcal{A}_{i}\theta_{i} = b_{i}
\end{equation}
where $\theta=[Q_s^\mathrm{T},(R_i)_s^\mathrm{T}]^\mathrm{T}$ and $b=(C_{e}(N))_{s}$ with $C_{e}(N) = [C_{e,0},C_{e,1}^\mathrm{T},\ldots,C_{e,N-1}^\mathrm{T}]^\mathrm{T}$.
$\mathcal{A}_{i}$ satisfies
\begin{equation}\label{A_LS_BCI}
\begin{aligned}
\mathcal{A}_{i} &= [D_{i},D_{i}(FK_{i}\otimes FK_{i})+(I_{n_{x}}\otimes \Gamma_{i})] \\
  D_{i} &= (H_{i}\otimes\mathcal{O}_{i})(I_{n_{x}^{2}}-\bar{F}_i\otimes\bar{F}_i)^{-1} \\
  \mathcal{O}_{i} &= \big[H_{i}^\mathrm{T},(H_{i}\bar{F}_i)^\mathrm{T},\ldots,(H_{i}\bar{F}_i^{N-1})^\mathrm{T}\big]^\mathrm{T} \\
  \Gamma_{i} &= [I_{n_{z}},-(H_{i}FK_{i})^\mathrm{T},\ldots,-(H\bar{F}_i^{N-2}FK_{i})^\mathrm{T}]^\mathrm{T}
  \end{aligned}
\end{equation}

The parameters $\theta$ is computed as the solution of semi-definite constrained least squares problem
\begin{equation}\label{ALS}
  \hat{\theta}_{i} = \argmin_{\theta_{i}}\|\mathcal{A}_{i}\theta_{i} - b_{i}\|_{2}^2 \quad s.t.,\quad Q,R_{i} \ge 0
\end{equation}
where the matrix inequalities $Q,R_{i}\ge 0$ can be handled by adding a logarithmic barrier function to the objective.
\cite{Rajamani2009} proves the uniqueness of the solution to the problem is guaranteed if and only if $\mathcal{A}$ has full column rank.
Furthermore, if $(F,H)$ is observable and $F$ is non-singular, the optimization in (\ref{ALS}) has a unique solution if and only if $\dim[null(D)]=0$.

It should be noted that when the dimension of the state $x$ is large and the window size of auto-covariance is small, the equation (\ref{ALS}) is easy to fall into overfitting problem.
To alleviate it, the $L_{2}$ regularization term is applied to (\ref{ALS}), then
\begin{equation}\label{SDP_ALS}
    \hat{\theta}_{i} = \argmin_{\theta_{i}}\|\mathcal{A}_{i}\theta_{i} - b_{i}\|_{2}^2 + \mu \|\theta_{i}\|_{2}^2
\end{equation}
where $\|\theta\|_{2}$ can be replaced by the trace of process noise covariance $tr(Q)$ for simplicity.
$\mu$ is the regularization term, and a good value of $\mu$ is such that $tr(Q)$ is small and any further decrease in value of $tr(Q)$ causes significant increase.
When the matrix inequality holds, $\hat{\theta}_{i}$ is estimated in the minimum mean-square error sense as
\begin{equation}\label{ALS_estimate}
  \hat{\theta}_{i} = (\mathcal{A}_{i}^\mathrm{T}\mathcal{A}_{i}+\mu I)^{-1}\mathcal{A}_{i}^\mathrm{T}\hat{b}_{i}=\mathcal{A}_{i}^{+}\hat{b}_{i}
\end{equation}
where $\hat{b}=(\hat{C}_{e}(N))_{s}$ is the unbiased estimate of the vector $b$ and computed as the empirical mean of the $i^{th}$ agent's auto-covariance innovations $\hat{C}_{e,i,l}$ is computed by using the ergodic property of the $L$-innovations from the given set of data
$\hat{C}_{e,i,l} = \frac{1}{\tau-l} \sum_{k=1}^{\tau-l}e_{i,k+l}e_{i,k}^\mathrm{T}$.

It is shown in \cite{Dunik2016} that the optimal estimator gain $K_{i}$ can be determined as $K_{i}^{\star} = \argmin_{K_{i}} f(\mathcal{J}(K_{i}))$ by minimizing the upper bound of the variance of the ALS estimate according to Isserlis' theorem, denoted as $\hat{\theta}$ $P_{\hat{\theta}}=cov[\hat{\theta}]=\mathbb{E}[(\theta-\hat{\theta})(\theta-\hat{\theta})^\mathrm{T}]=\mathcal{A}^{+}cov[\hat{b}]\mathcal{A^{+}}^\mathrm{T}$.
where $f(\cdot)$ is a suitable function, e.g., the trace.
$\mathcal{J}(K_{i})$ is the known criterion of $K_{i}$ defined in \cite{Dunik2016}.
\section{Distributed ALS method}
\label{sec:distributeALS}
The ALS estimator of the noise covariance matrices is proven to be unbiased and converging asymptotically to the true values with increasing number of data $\tau$ in \cite{Odelson2006}.
But in sensor network, each agent has a limited storage capacity and suffers from the lack of innovations, i.e., $\tau$ is small.
Besides, with the increase of $\tau$, the computation and stoage burden of each agent will become heavier.
Therefore, it is necessary to reformulate the distributed ALS method to balance the tradeoff between the state estimation accuracy and the computation capacity.
One effective approach is that each agent enlarges the number of input data $\tau$ by receiving the auto-covariance from its neighbors $j\in N_{i}\cup i$.
The equation (\ref{SDP_ALS}) turns into a joint cost function, as shown below.
\begin{equation}\label{alg_opt_distest_Q}
  \hat{\theta}_{i} =\argmin_{\theta_{i}} \sum_{j\in N_{i}\cup i}\Big(\|\mathcal{A}_{j}\theta_{j}-b_{j}\|_{2}^2 + \mu \|\theta_{j}\|_{2}^2
\end{equation}
The empirical mean of the auto-covariance innovations term $\hat{C}_{e,i,l}$ is reformulated as the mean of the neighbors and itself
\begin{equation}\label{Cfusion}
  \hat{C}_{e,i,l} = \frac{1}{\tau-l} \frac{1}{M_{i}+1}\sum_{j\in N_{i}\cup i}\sum_{k=1}^{\tau-l}e_{j,k+l}e_{j,k}^\mathrm{T}
\end{equation}

\begin{figure}[H]
  \centering
  \includegraphics[width=\hsize]{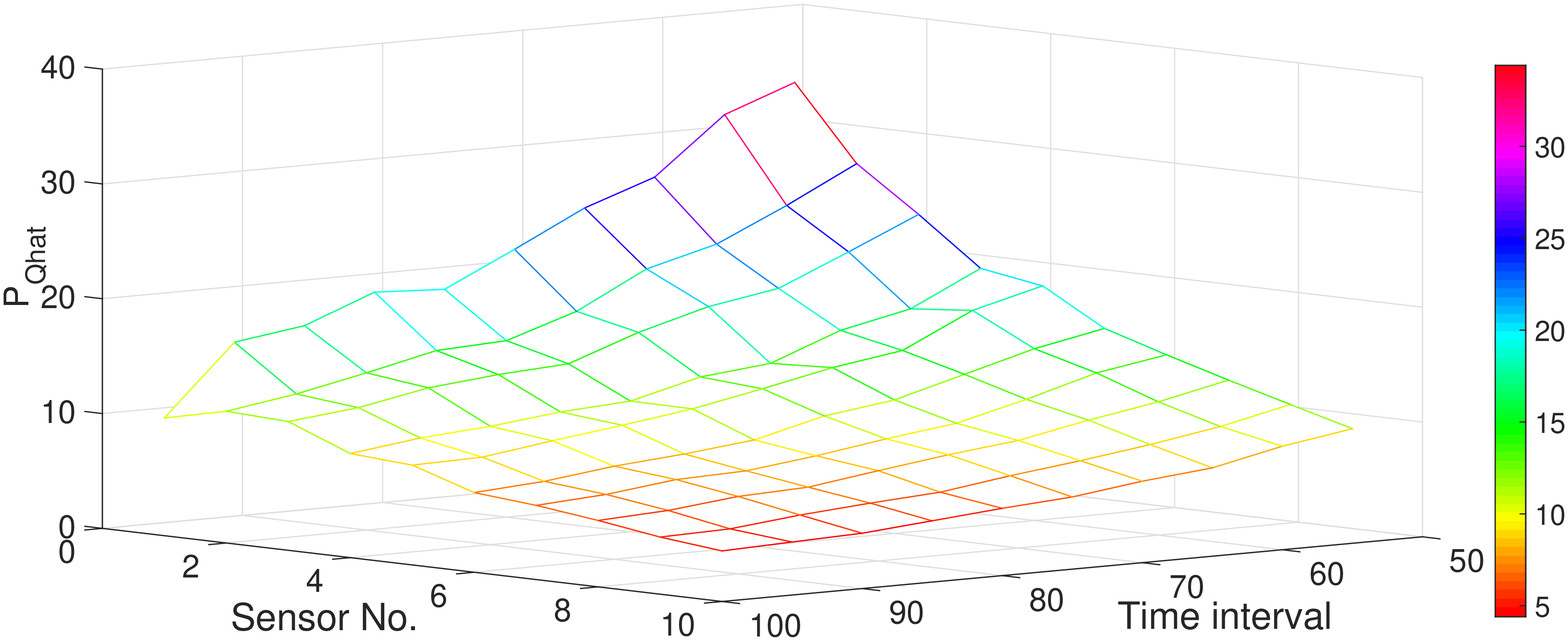}
  \caption{Relationship between $var(\hat{Q}_{i})$ and $(\mathcal{N}_{i},\tau)$}\label{Qvar}
\end{figure}

We do experiments on a linear time-invariant system with $10$ sensors to be deployed in a fully connected network.
The variance of the ALS-estimates $var(\hat{\theta}_{i})=var(\hat{Q}_{i})$ denotes the state estimation accuracy.
$F=-0.8$, $H_{i}=1, i=1,\ldots,10$, and with true but unknown noise variances $Q=8$ and $R=[1,2,\ldots,10]$.
The time intervals are set as $\tau=[55,60,\ldots,100]$.
The system is simulated for $10^4$ steps.
The relationship between $var(\hat{Q}_{i})$ and the number of sensors $\mathcal{N}_{i}$ to be fused plus the different values of $\tau$ is shown in Fig. \ref{Qvar}.

As is indicated from Fig. \ref{Qvar}, the variance of ALS estimate $var(\hat{Q})$ decreases with the increase of the number of sensors $\mathcal{N}_{i}$ and the number of innovations.
$var(\hat{Q})$ decreases from 34.4 to 10.25 as the number of innovations increases from 55 to 100 when the number of sensors is 1,
and $var(\hat{Q})$ decreases from 10.25 to 4.97 as the number of fused sensors increases from 1 to 10 when the number of innovations is 100.
Therefore, it is possible for each agent to reduce the variance of the ALS estimate by receiving the innovations from its neighbors instead of increasing the number of innovations.

The empirical mean of the auto-covariance innovations term $\hat{C}_{e,i,l}$ in (\ref{Cfusion}) is only the fusion of $\hat{b}_{i}$.
To derive the optimal fused noise covariance estimate denoted as $\hat{\theta}_{F}^{\star}$, the fused residual $\varepsilon_{k,F}$ and residual covariance $P_{\varepsilon,k,F}$ is computed based on the batch covariance intersection (BCI) method:
\begin{equation}\label{BCI_P}
  P_{\varepsilon,k,F}^{-1} = \sum_{j\in N_{i}\cup i}w_{j}P_{\varepsilon,k,j}^{-1}
\end{equation}
\begin{equation}\label{BCI_x}
\begin{aligned}
  P_{\varepsilon,k,F}^{-1}\varepsilon_{k,F} &= \sum_{j\in N_{i}\cup i}w_{j}P_{\varepsilon,k,j}^{-1}\varepsilon_{k,j}\\
  \sum_{i=1}^{N} w_{i} &= 1, w_{i}\in [0,1],i=1,2,\ldots,N
  \end{aligned}
\end{equation}
where the weights $w_{i}$ can be determined by using some sub-optimal methods such as minimizing the trace of fused residual covariance $P_{\varepsilon,k,F}^{-1}$ in \cite{Franken2005}.
\begin{equation}\label{BCI_w}
  w_{i} = \frac{1/tr(P_{\varepsilon,k,i})}{\sum\limits_{j\in N_{i}\cup i}1/tr(P_{\varepsilon,k,j})}
\end{equation}
\begin{equation}\label{CeN_BCI}
P_{\varepsilon,k,i} = \Big(\sum_{j\in N_{i}\cup i}\frac{1/tr(P_{\varepsilon,k,i})}{\sum\limits_{j\in N_{i}\cup i}1/tr(P_{\varepsilon,k,j})}P_{\varepsilon,k,j}^{-1}\Big)^{-1}
\end{equation}

Denote the matrices $\mathcal{A}_{F}$ and $\hat{b}_F$ as $\mathcal{A}_{F} = \bigoplus_{i=1}^{M_{i}}\mathcal{A}_{i}$ and $\hat{b}_F = \hat{b}_{i} \otimes I_{M_{i}}$, where $\hat{b}_{i}=[\hat{C}_{e,i,0}^\mathrm{T},\hat{C}_{e,i,1}^\mathrm{T},\ldots,\hat{C}_{e,i,N-1}^\mathrm{T}]^\mathrm{T}$.

Then the solution to the problem in (\ref{alg_opt_distest_Q}) can be solved by solving the regularized LS problem.
\begin{equation}\label{Ab_LS_combine}
  \hat{\theta}_F^{\star} = \argmin_{\theta_F}\|\mathcal{A}_{F}\theta_F-\hat{b}_F\|_{2}^2 + \mu \|\theta_F\|_{2}^2
\end{equation}
where $\hat{\theta}_F = \hat{\theta}_{i} \otimes I_{M_{i}}$.
The problem is solved as
\begin{equation}\label{theta_combine}
  \hat{\theta}_F^{\star} = (\mathcal{A}_F^\mathrm{T}\mathcal{A}_F+\mu I)^{-1}\mathcal{A}_F^\mathrm{T}\hat{b}_F=\mathcal{A}_F^{+}\hat{b}_F
\end{equation}

Then the ALS method combined with the BCI algorithm is summarized in Alg. \ref{alg_als_bci}.
\begin{algorithm}[H]
\caption{\small Solving problem (\ref{alg_opt_distest_Q}) by D-ALS algorithm}
\label{alg_als_bci}
\textbf{Input:} $\mu=0.01$, $\nu=5\times10^{-3}$, $\tau=100$, $N_{sim}=10^3$. \\ 
\textbf{Initialize:} $k=0$, $\hat{x}_{0}$, $Q_0$, $z_{i,1:N_{sim}}$, $R_{0,i}$, $P_{\varepsilon,0}$ and $K_{0}$, $i=1,\cdots,M$.\\
\textbf{Output:} $\hat{Q}^{\star}$. \\
\textbf{while} in loop and $\hat{Q}_{k+1}-\hat{Q}_{k+1}> \nu$ \textbf{do}
\begin{enumerate}
  \item Update $\hat{x}_{i,k+1}$, $\varepsilon_{i,k+1}$, $P_{\varepsilon,i}$ and $K_{i}$ in (\ref{eps_i}) to (\ref{Peps_i}), and then calculate the fused residual and its covariance $P_{\varepsilon,k,F}$ by BCI method in (\ref{BCI_P}) and (\ref{BCI_x}). Update the matrix $\mathcal{A}_{i}$, $\hat{b}_{i}$, $\mathcal{A}_{F}$ and $\hat{b}_{F}$ in .
  \item Update the global optimal noise covariance $\hat{\theta}^{\star}$ in (\ref{Ab_LS_combine}) to (\ref{theta_combine}), and set $\hat{Q}_{k+1}=\hat{\theta}_{1}^{\star}$.
\end{enumerate}
\textbf{end while}
\label{alg_lirnn}
\end{algorithm}
\begin{remark}
It is easily derived that
the augmented matrix $\mathcal{A}_F$ and the permutation matrix $\hat{b}_{F}$ for $i^{th}$ sensor has dimensions of $M_iNn_{x}\times (n_{x}+M_{i}n_{z})$ and $M_iNn_{x}(n_{x}+M_{i}n_{z})\times 1$ respectively.
The computation complexity of the D-ALS algorithm is $O(M_iN^2n_{x})$.
Expanding the number of the auto-covariance of innovations would increase the computation time and even lead to the intractable computation.
Therefore, the number of sensors and the window size should be made from a tradeoff between the accuracy and the computation burden.
\end{remark}
%
%
Then the variance of the fused noise covariance matrix $\hat{\theta}_{F}$ is lower than the variance of each agent's noise covariance matrix $\hat{\theta}_{i}$, as is proved in Theorem \ref{theta_proof}.
\begin{theorem}\label{theta_proof}(Accuracy Analysis of D-ALS algorithm)
The relations between local and fused residuals covariance $P_{\varepsilon,k,i}$, $P_{\varepsilon,k,0}$, $\bar{P}_{\varepsilon,k,F}$ and $P_{\varepsilon,k,F}$ are shown as follows.
  \begin{equation}\label{P_{0}}
    tr(P_{\varepsilon,k,0})\le tr(\bar{P}_{\varepsilon,k,F}) \le tr(P_{\varepsilon,k,F}) \le tr(P_{\varepsilon,k,i})
  \end{equation}

Then the relations between the fused noise covariance matrix $\hat{\theta}_{F}$ and the noise covariance matrix for each agent are shown below.
\begin{equation}
  var(\hat{\theta}_{F}) \le var(\hat{\theta}_{i}), \quad i=1,\ldots,M_{i}
\end{equation}
\end{theorem}
\textbf{Proof}:
Using the unbiasedness of $\hat{\varepsilon}_{k,i}$ for each agent $i$, it can be derived that $\hat{\varepsilon}_{F,k}$ is a linear unbiased estimate.
Since $\hat{\varepsilon}_{k,0}$ is the best linear unbiased estimate, then $P_{\varepsilon,k,0}\le \bar{P}_{\varepsilon,k,F}$ holds.
The inequality $\bar{P}_{\varepsilon,k,F} \le P_{\varepsilon,k,F}$ is proved as the consistency property in \cite{Julier1997}.
Because the operator of trace is monotonically increasing function, then the inequalities $tr(P_{\varepsilon,k,0})\le tr(\bar{P}_{\varepsilon,k,F})$ and $tr(\bar{P}_{\varepsilon,k,F}) \le tr(P_{\varepsilon,k,F})$ hold.
When the parameters are set to $w_{i}=1$ and $w_{j}=1,j\neq i$, $tr(P_{\varepsilon,k,F})=tr(P_{\varepsilon,k,i})$.
Because the parameter $w$ is determined by minimizing the trace of $P_{\varepsilon,k,F}$, as shown in (\ref{w_proof}), it is easily derived that $P_{\varepsilon,k,F}\le P_{\varepsilon,k,i}$.
\begin{equation}\label{w_proof}
  w = \argmin_{w} tr\Big[(\sum_{i=1}^{M_{i}}w_{i}P_{\varepsilon,k,i}^{-1})^{-1}\Big]
\end{equation}

Then it can easily be derived that $\mathcal{A}_{F}$ is the best estimate of $\mathcal{A}$.
As shown in (\ref{Cfusion}), the number of data to compute $\hat{b}_F$ is larger than that to compute $\hat{b}_i$.
As is proven in \cite{Odelson2006}, The ALS estimate of the noise covariance matrix converges asymptotically to the true values with increasing number of data.
Therefore, the variance of $\hat{\theta}_{F}$ is smaller than the variance of $\hat{\theta}_{i}$.
Q.E.D.
\section{Simulation Results}\label{sec:simulation}
\subsection{Static sensor networks}\label{SSN}
Consider fully connected sensor network
The linear time-invariant system is modeled as $x_{k+1} = 0.8x_{k} + w_{k}$ measured by $3$ sensors modeled as $y_{i,k+1} = H_{i}x_{k} + v_{i,k}$
with $H_{1} = [1,0]$, $H_{2} = I_{2}$, $H_{3} = [1,0]$.
$w_{k}$ and $v_{i,k}$ are zero-mean Gaussian noise with unknown covariance $Q$ and $R_{i}$, where the real values are set as $Q=4$, $R_{1}=0.81$, $R_{2}=diag(4,0.64)$, and $R_{3}=2.25$.
Here, to guarantee the unbiasedness, we ran $N_{s}=10^4$ Monte Carlo simulations for each simulation data set.
The estimation performance for the $i_{th}$ sensor is measured as the mean square error (MSE):
\begin{equation}\label{defMSE}
    MSE_{i,k} = \frac{1}{N_{s}}\sum_{t=k}^{k+N_{s}}(\varepsilon_{i,t} - \hat{\varepsilon}_{i,t})^{2}, \quad i=1,2,,3,M_{i}
\end{equation}

\begin{figure}
  \centering
  \includegraphics[width=\hsize]{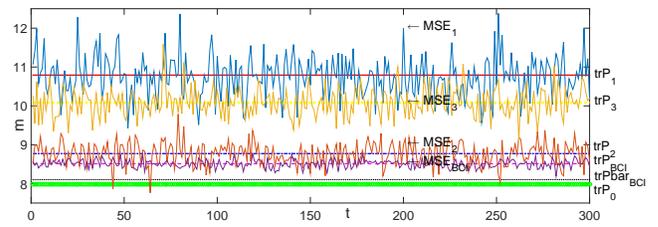}
  \caption{Comparisons of $MSE$ curve and the trace of residual covariance between fused and single sensor.}\label{Comp_Peps}
\end{figure}

Using the steady-state Kalman filter, the gain $K_{i}$ is computed as
\begin{equation}\label{K_sim}
  K_{1} = [1.68,0.81], K_{2} = \left[
                                       \begin{array}{cc}
                                         0.36 & 1.22 \\
                                         0.06 & 0.84 \\
                                       \end{array}
                                     \right], K_{3} = [1.40,0.61]
\end{equation}
According to the equations (\ref{BCI_P}), the steady state covariances of each agent and the cross covariance matrices are
\begin{equation}\label{P_sim}
\begin{aligned}
  &P_{\varepsilon,1} = \left[
                                       \begin{array}{cc}
                                         5.05 & 4.94 \\
                                         4.94 & 5.73 \\
                                       \end{array}
                                     \right], \quad P_{\varepsilon,2} = \left[
                                       \begin{array}{cc}
                                         3.73 & 2.78 \\
                                         2.78 & 5.05 \\
                                       \end{array}
                                     \right] \\
  &P_{\varepsilon,3} = \left[
                                       \begin{array}{cc}
                                         5.79 & 3.65 \\
                                         3.65 & 4.29 \\
                                       \end{array}
                                     \right], \quad
  P_{\varepsilon,1,2} = \left[
                                       \begin{array}{cc}
                                         1.30 & -0.22 \\
                                         -0.22 & 0.38 \\
                                       \end{array}
                                     \right] \\
  &P_{\varepsilon,2,3} = \left[
                                       \begin{array}{cc}
                                         0.75 &  0.21 \\
                                         0.21 & 0.45 \\
                                       \end{array}
                                     \right], \quad
  P_{\varepsilon,1,3} = \left[
                                       \begin{array}{cc}
                                         0.69 & 1.21 \\
                                         1.21 & 4.19 \\
                                       \end{array}
                                     \right] \\
\end{aligned}
\end{equation}
The comparisons between MSE and the trace of $P_{\varepsilon,k,i}$, $P_{\varepsilon,k,0}$, $\bar{P}_{\varepsilon,k,F}$ and $P_{\varepsilon,k,F}$ are shown in Fig. \ref{Comp_Peps} and Table \ref{Comp_P}.
\begin{table}
\caption{Comparisons of the trace of residual covariance between fused and single sensor.}
\centering
\begin{tabular}{cccccc}
\toprule
$trP_{1}$ & $trP_{2}$ & $trP_{3}$ & $trP_{F}$ & $tr\bar{P}_{F}$ & $tr\bar{P}_{0}$ \\
\midrule
10.791 & 10.083 & 8.771 & 8.512 & 8.102 & 7.985\\
\bottomrule
\end{tabular}
\label{Comp_P}
\end{table}
As indicated from the figure, the true accuracy of the BCI fused residual covariance is similar to the linear optimal residual covariance, because $tr(\bar{P}_{\varepsilon,k,F})=8.102$ is near to $tr(\bar{P}_{\varepsilon,k,0})=7.985$.
Besides, the variance of the fused estimate $tr\bar{P}_{F}$ is lower than others, illustrating that the proposed D-ALS algorithm outperforms than the ALS method.
\subsection{Mobile sensor networks}
Cooperative target tracking in mobile sensor networks (MSNs) is an important task in many applications, e.g., the unmanned aerial vehicle (UAV).
Compared with the target tracking case in static sensor networks (SSNs) in section \ref{SSN}, each agent node in this case is mobile and versatile, and is required to be deployed in any scenario with rapid topology changes.
Therefore, the ALS-BCI algorithm is also applied to the target tracking in MSNs to show its efficiency.

Consider 10 sensor nodes tracking the maneuver target in a $110m\times90m$ square.
The target is driven by a turning rate model:
\begin{equation}\label{MSN_target}
  X_{k+1} = \left[
                  \begin{array}{cccc}
                    1 & \frac{\sin(\eta T_{s})}{\eta} & 0 & -\frac{1-\cos(\eta T_{s})}{\eta} \\
                    0 & \cos(\eta T_{s}) & 0 & -\sin(\eta T_{s}) \\
                    0 & \frac{1-\cos(\eta T_{s})}{\eta} & 1 & \frac{\sin(\eta T_{s})}{\eta} \\
                    0 & \sin(\eta T_{s}) & 0 & \cos(\eta T_{s}) \\
                  \end{array}
                \right]X_{k}+Gw_{k}
\end{equation}
where $X_{k}=[x_{k},\dot{x}_{k},y_{k},\dot{y}_{k}]^\mathrm T$ is the states to be estimated at time $k$.
The states includes the position $[x_{k},y_{k}]$ and the velocity $[\dot{x}_{k},\dot{y}_{k}]$, and the initial values are $[10m,2m/s,100m,2m/s]$.
$\eta$ is the turn rate and is set to $\frac{\pi}{60}rad/s$.
$w_{k}$ is Gaussian white noise with covariance matrix $Q=diag[Q_x,Q_{\dot{x}},Q_y,Q_{\dot{y}}]$, where $Q_x,Q_y,Q_{\dot{x}},Q_{\dot{y}}$ are unknown scalar variable to be estimated and the real ones are set to $10m^2$, $0$, $10m^2$ and $0$  respectively.
$T_{s}$ is the sampling time and is set to $1s$.
$G=\left[
     \begin{array}{cccc}
       \frac{1}{2}T_{s}^2 & T_{s} & 0 & 0 \\
       0 & 0 & \frac{1}{2}T_{s}^2 & T_{s}\\
     \end{array}
   \right]^\mathrm T
$.

The measurements of each agent is given by:
\begin{equation}\label{WSN_M}
  Z_{i,k}=\left[
                \begin{array}{cccc}
                  1 & 0 & 0 & 0 \\
                  0 & 0 & 1 & 0 \\
                \end{array}
              \right]X_{k} + v_{i,k} \quad i=1,\ldots,M
\end{equation}
where $Z_{i,k}=[zx_{k},zy_{k}]^T$ is the measurement of the position of the target.
$v_{i,k}$ is Gaussian white noise with unknown covariance matrix $R_{i}=diag[R_{x,i},R_{y,i}]$.
$R_{x,i}$ and $R_{y,i}$ are unknown parameter to be estimated and the real ones are set to $2m^2$ and $2m^2$ respectively.
The motion of each agent is described by the following kinematic equations:
\begin{equation}\label{model}
\begin{aligned}
qx_{i}(k) &= qx_{i}(k-1) + qv_{i}*T_{s}*\cos \theta_{k} \\
qy_{i}(k) &= qy_{i}(k-1) + qv_{i}*T_{s}*\sin \theta_{k}
\end{aligned}
\end{equation}
where $(qx_{i}(k),qy_{i}(k))$ is the position of the $i^{th}$ sensor at time $k$.
$\theta_{k}=\arctan \frac{z_{i,k}-qy_{i}(k-1)}{z_{i,k}-qx_{i}(k-1)}$ is the measurement of angular position of the $i^{th}$ sensor towards the target.
$qv_{i}$ is the constant speed of the $i^{th}$ sensor and is set to $0.5m/s$.
The initial position of the sensors are set to $(18,27)$, $(31,43)$, $(62,41)$, $(86,33)$, $(15,45)$, $(13,98)$, $(38,105)$, $(60,99)$, $(89,93)$, $(110,106)]$, and the unit is meter.
The communication ranges and the sensing ranges of the sensors are all set to $r_{c}=45m$ and $r_{s}=60m$ respectively.

\begin{figure}
  \centering
  \includegraphics[width=\hsize]{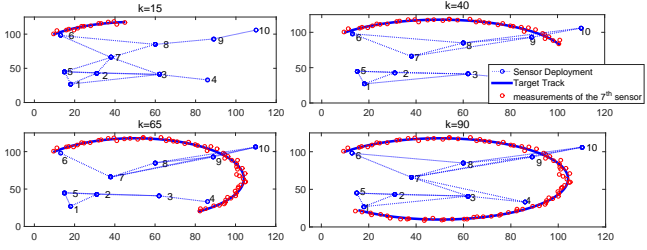}
  \caption{The diagram of the target tracking in the time-varying sensor deployment.}\label{WSN_deploy}
\end{figure}

\begin{figure}
  \centering
  \includegraphics[width=\hsize]{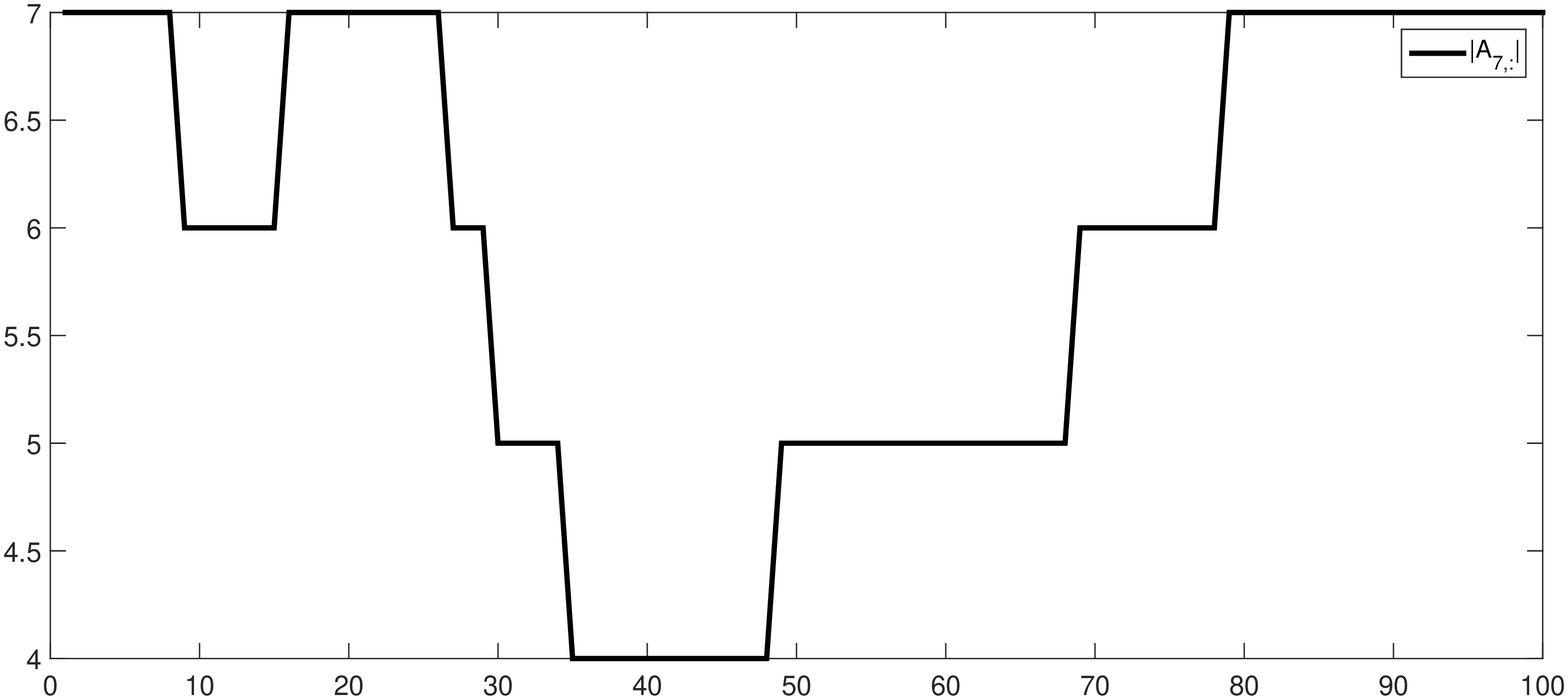}
  \caption{The connectivity of the $7^{th}$ sensor.}\label{WSN_connect}
\end{figure}

100 Monte Carlo simulations are run on the simulated model.
For comparison we ran ALS and ALS-BCI on the same data sets.
The initial parameters of the algorithm \ref{alg_opt_distest_Q} are set to $Q_{A}=diag[5,0,5,0]$ and $R_{A}=diag[1,1]$.
The diagram of the target tracking in the time-varying sensor deployment is shown in Fig. \ref{WSN_deploy}.
The connectivity of the $i^{th}$ sensor is denoted as the sum of the adjacent matrix $|A_{ij}|$, as shown in Fig. \ref{WSN_connect}.

As is shown from the above figures, the $7^{th}$ sensor only measures the target in the sensing range, and the communication topology is varying due to the mobility of the sensors.

\begin{figure}[H]
  \centering
  \includegraphics[width=\hsize]{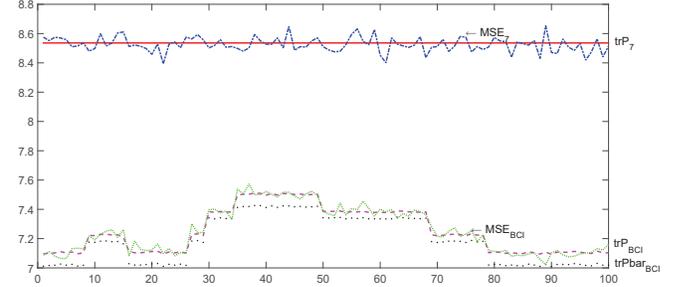}
  \caption{Comparisons of $MSE$ curve and the residual covariance between fused and $7^{th}$ sensor.}\label{Comp_Peps_WSN}
\end{figure}

The comparisons between MSE and the trace of the $7^{th}$ sensor's residual covariance $P_{7}$ and its fused residual covariance $P_{BCI}$ and $\bar{P}_{BCI}$ along the x-axis are shown in Fig. \ref{Comp_Peps_WSN} and Table \ref{Comp_P_WSN}.
\begin{table}[H]
\caption{Comparisons of the residual covariance between fused and $7^{th}$ sensor.}
\centering
\begin{tabular}{c|cccc}
\toprule
\multirow{2}{*}{$trP_{7}$} & \multicolumn{2}{c}{$trP_{BCI}$} & \multicolumn{2}{c}{$tr\bar{P}_{BCI}$}\\
 & k=15 & k=43  & k=15 & k=43\\
\midrule
8.536 & 7.223 & 7.502 & 7.185 & 7.411\\
\bottomrule
\end{tabular}
\label{Comp_P_WSN}
\end{table}

As is indicated from the figure and the table, the fused residual covariance $trP_{BCI}$ and the MSE values $MSE_{BCI}$ of $7^{th}$ sensor is less than that of its own $trP_{7}$ and $MSE_{7}$.
Besides, it is noticed that $trP_{BCI}$ and $tr\bar{P}_{BCI}$ decreases with the increase of the number of neighbors.
\section{CONCLUSIONS AND FUTURE WORKS}

\subsection{Conclusions}
\label{sec:conclusions}
This paper proposes the distributed auto-covariance least squares algorithm based on BCI to solve the distributed estimation problem with unknown noise covariance over networks.
The efficiency of the algorithm is proven because the fused error covariances converges to the true values faster and the variance of the ALS estimate is smaller.
The numerical results are illustrated to show the performance of the algorithm.
\subsection{Future Works}
In real-time applications, the latency and limited power are the main problems in wireless sensor network. Since the ALS-BCI algorithm still needs some time to compute the matrix $\mathcal{A}$, it is necessary to compute the a-priori estimate of the lower bound of the variance of the fused ALS estimate $\hat{\theta}_{F}$.
Future work is needed to derive the lower bound of the fused noise covariance estimate $\hat{\theta}_{F}$.



\bibliographystyle{IEEEtran}
\bibliography{myCDCref}

\begin{thebibliography}{10}
\providecommand{\url}[1]{#1}
\csname url@samestyle\endcsname
\providecommand{\newblock}{\relax}
\providecommand{\bibinfo}[2]{#2}
\providecommand{\BIBentrySTDinterwordspacing}{\spaceskip=0pt\relax}
\providecommand{\BIBentryALTinterwordstretchfactor}{4}
\providecommand{\BIBentryALTinterwordspacing}{\spaceskip=\fontdimen2\font plus
\BIBentryALTinterwordstretchfactor\fontdimen3\font minus
  \fontdimen4\font\relax}
\providecommand{\BIBforeignlanguage}[2]{{%
\expandafter\ifx\csname l@#1\endcsname\relax
\typeout{** WARNING: IEEEtran.bst: No hyphenation pattern has been}%
\typeout{** loaded for the language `#1'. Using the pattern for}%
\typeout{** the default language instead.}%
\else
\language=\csname l@#1\endcsname
\fi
#2}}
\providecommand{\BIBdecl}{\relax}
\BIBdecl

\bibitem{He2019}
S.~He, H.~Shin, S.~Xu \emph{et~al.}, ``Distributed estimation over a low-cost
  sensor network: A review of state-of-the-art,'' \emph{Information Fusion},
  vol.~54, pp. 21 -- 43, 2020.

\bibitem{Cattivelli2010}
F.~S. Cattivelli and A.~H. Sayed, ``Diffusion strategies for distributed
  {K}alman filtering and smoothing,'' \emph{IEEE Trans. on Automat. Cont.},
  vol.~55, no.~9, pp. 2069--2084, 2010.

\bibitem{Chong2001}
C.~Chong, B.~Allen, I.~Hamilton \emph{et~al.}, ``Convex combination and
  covariance intersection algorithms in distributed fusion,'' in \emph{Proc.
  Fusion'01}, 2001.

\bibitem{Sun2016}
T.~Sun, M.~Xin, and B.~Jia, ``Distributed estimation in general directed sensor
  networks based on batch covariance intersection,'' in \emph{American Control
  Conference, 2016.}, 2016, pp. 5492--5497.

\bibitem{TerryMoore2003}
TerryMoore and MartinSmith, ``Adaptive {K}alman filtering for low-cost
  {INS/GPS},'' \emph{Journal of Navigation}, vol.~56, no.~1, pp. 143--152,
  2003.

\bibitem{Deng2017}
F.~{Deng}, S.~{Guan}, X.~{Yue} \emph{et~al.}, ``Energy-based sound source
  localization with low power consumption in wireless sensor networks,''
  \emph{IEEE Trans. on Industrial Electronics}, vol.~64, no.~6, pp. 4894--4902,
  2017.

\bibitem{Gertler1988}
J.~J. Gertler, ``Survey of model-based failure detection and isolation in
  complex plants,'' \emph{IEEE Control Systems Magazine}, vol.~8, no.~6, pp.
  3--11, 1988.

\bibitem{Mehra1970}
R.~K. Mehra, ``On the identification of variances and adaptive {K}alman
  filtering,'' \emph{IEEE Trans. on Automat. Cont.}, vol.~15, no.~2, pp.
  175--184, 1970.

\bibitem{Belanger1974}
P.~R. B\'{e}langer, ``Estimation of noise covariance matrices for a linear
  time-varying stochastic process,'' \emph{Automatica}, vol.~10, no.~3, pp.
  267--275, 1974.

\bibitem{Odelson2006}
B.~J. Odelson, M.~R. Rajamani, and J.~B. Rawlings, ``A new autocovariance
  least-squares method for estimating noise covariances,'' \emph{Automatica},
  vol.~42, no.~2, pp. 303--308, 2006.

\bibitem{Akesson2008}
B.~M. \r{A}kesson, J.~B. Jørgensen, N.~K. Poulsen \emph{et~al.}, ``A
  generalized autocovariance least-squares method for {K}alman filter tuning,''
  \emph{Journal of Process Control}, vol.~18, no. 7-8, pp. 769--779, 2008.

\bibitem{Rajamani2009}
M.~R. Rajamani and J.~B. Rawlings, ``Estimation of the disturbance structure
  from data using semidefinite programming and optimal weighting,''
  \emph{Automatica}, vol.~45, no.~1, pp. págs. 142--148, 2009.

\bibitem{Dunik2009}
J.~Dunik, ``Methods for estimating state and measurement noise covariance
  matrices: Aspects and comparison,'' in \emph{System Identification}, 2009,
  pp. 372--377.

\bibitem{Abdel2010}
M.~F. Abdel-Hafez, ``The autocovariance least-squares technique for {GPS}
  measurement noise estimation,'' \emph{IEEE Trans. on Vehicular Technology},
  vol.~59, no.~2, pp. 574--588, 2010.

\bibitem{Lima2013}
F.~V. Lima, M.~R. Rajamani, T.~A. Soderstrom \emph{et~al.}, ``Covariance and
  state estimation of weakly observable systems: Application to polymerization
  processes,'' \emph{IEEE Trans. on Control Systems Technology}, vol.~21,
  no.~4, pp. 1249--1257, 2013.

\bibitem{Dunik2016}
J.~Dunik, O.~Straka, and M.~Simandl, ``On autocovariance least-squares method
  for noise covariance matrices estimation,'' \emph{IEEE Trans. on Automat.
  Cont.}, vol.~62, no.~2, pp. 967--972, 2017.

\bibitem{Sarkka2009}
S.~S\"{a}rkk\"{a} and A.~Nummenmaa, ``Recursive noise adaptive {K}alman
  filtering by variational {Bayesian} approximations,'' \emph{IEEE Trans. on
  Automat. Cont.}, vol.~54, no.~3, pp. 596--600, 2009.

\bibitem{Li2019}
J.~{Li}, F.~{Deng}, and J.~{Chen}, ``A fast distributed variational bayesian
  filtering for multisensor ltv system with non-gaussian noise,'' \emph{IEEE
  Trans. on Cybernetics}, vol.~49, no.~7, pp. 2431--2443, July 2019.

\bibitem{Kashyap1970}
R.~Kashyap, ``Maximum likelihood identification of stochastic linear systems,''
  \emph{IEEE Trans. on Automat. Cont.}, vol.~15, no.~1, pp. 25--34, Feb 1970.

\bibitem{Chen2017}
J.~{Chen}, J.~{Li}, S.~{Yang} \emph{et~al.}, ``Weighted optimization-based
  distributed {K}alman filter for nonlinear target tracking in collaborative
  sensor networks,'' \emph{IEEE Trans. on Cybernetics}, vol.~47, no.~11, pp.
  3892--3905, 2017.

\bibitem{Myers1976}
K.~Myers and B.~Tapley, ``Adaptive sequential estimation with unknown noise
  statistics,'' \emph{IEEE Trans. on Automat. Cont.}, vol.~21, no.~4, pp.
  520--523, Aug 1976.

\bibitem{Verdu1984}
S.~Verdu and H.~Poor, ``Minimax linear observers and regulators for stochastic
  systems with uncertain second-order statistics,'' \emph{IEEE Trans. on
  Automat. Cont.}, vol.~29, no.~6, pp. 499--511, Jun 1984.

\bibitem{Overschee1996}
P.~V. Overschee and B.~D. Moor, ``Subspace identification for linear systems:
  Theory–implementation–applications,'' \emph{Kluwer Academic Publishers},
  pp. 57--93, 1996.

\bibitem{Carew1974}
B.~Carew and P.~Belanger, ``Identification of optimum filter steady-state gain
  for systems with unknown noise covariances,'' \emph{IEEE Trans. on Automat.
  Cont.}, vol.~18, no.~6, pp. 582--587, 1974.

\bibitem{Wiberg2000}
D.~M. Wiberg, T.~D. Powell, and D.~Ljungquist, ``An online parameter estimator
  for quick convergence and time-varying linear systems,'' \emph{IEEE Trans. on
  Automat. Cont.}, vol.~45, no.~10, pp. 1854--1863, 2000.

\bibitem{Basar1972}
T.~Basar and M.~Mintz, ``Minimax terminal state estimation for linear plants
  with unknown forcing functions,'' \emph{International Journal of Control},
  vol.~16, no.~16, pp. 49--69, 1972.

\bibitem{Tugnait1979}
J.~Tugnait and A.~Haddad, ``On state estimation for linear discrete-time
  systems with unknown noise covariances,'' \emph{IEEE Trans. on Automat.
  Cont.}, vol.~24, no.~2, pp. 337--340, 1979.

\bibitem{Simon2006}
D.~Simon, \emph{Optimal State Estimation: {K}alman, H Infinity, and Nonlinear
  Approaches}.\hskip 1em plus 0.5em minus 0.4em\relax Wiley-Interscience, 2006.

\bibitem{Franken2005}
D.~Franken and A.~Hupper, ``Improved fast covariance intersection for
  distributed data fusion,'' in \emph{International Conference on Information
  Fusion}, 2005.

\bibitem{Julier1997}
S.~J. Julier and J.~K. Uhlmann, ``A non-divergent estimation algorithm in the
  presence of unknown correlations,'' in \emph{American Control Conference,
  1997. Proceedings of the}, 1997, pp. 2369--2373 vol.4.

\end{thebibliography}

%
%
%

\end{document}